\documentstyle[aps]{revtex}   

\begin{document}
\draft
\title{ Einstein Metrics Adapted to Contact Structures on 3-Manifolds }
\author{Brendan S. Guilfoyle \footnote{Email: brendan.guilfoyle@ittralee.ie} }
\address{Mathematics Department, Institute of Technology Tralee, Tralee, Co. Kerry, Ireland.}
\date{\today}
\maketitle
\begin{abstract}
The Newman-Penrose-Perjes formalism is applied to smooth contact structures on riemannian 3-manifolds.  In particular it is shown that a contact 3-manifold admits an adapted riemannian metric if and only if it admits a metric with a divergence-free, constantly twisting, geodesic congruence. The shear of this congruence is identified with the torsion of the associated pseudohermitian structure, while the Tanaka-Webster curvature is identified with certain derivatives of the spin coefficients.  The particular case where the associated riemannian metric is Einstein is studied in detail.  It is found that the torsion is constant and the field equations are completely solved locally.  Hyperbolic space forms are shown not to have adapted contact structures, even locally, while contact structures adapted to a flat or elliptic space form are contact isometric to the standard one.
\end{abstract}

\newtheorem{thm}{Theorem}[section]   
\newtheorem{cor}[thm]{Corollary}   
\newtheorem{lem}[thm]{Lemma}   
\newtheorem{prop}[thm]{Proposition}   
\newtheorem{defn}{Definition}[section]      
\newtheorem{rem}{Remark}[section]   
\newcommand{\nc}{\newcommand}
\nc{\bee}{\begin{equation}}
\nc{\ee}{\end{equation}}
\nc{\beer}{\begin{eqnarray}}
\nc{\eer}{\end{eqnarray}}
\nc{\dv}{\textstyle\frac}
\nc{\sdv}{\scriptstyle\frac}
\nc{\hf}{\dv{1}{2}}
\nc{\shf}{\sdv{1}{2}}
\nc{\nn}{\nonumber}

\section{Introduction}

Over the last 30 years there has been increasing interest in contact structures from a number of perspectives.  In particular, contact structures have been found to be very useful for investigating the topology \cite{elias1} \cite{elias2} \cite{kanda} \cite{wein}  and geometry \cite{blair} \cite{geiges} \cite{gag} \cite{gold} \cite{kam} \cite{rud} \cite{tanno} of 3-manifolds.

A contact structure $\xi$ on a 3-manifold $M$ is maximally nonintegrable distribution of 2-planes.  Such a structure can always be defined by the vanishing of a (non-unique) 1-form $\alpha$, called the contact form.  On a contact 3-manifold Chern and Hamilton \cite{cah} introduced the definition of an {\it adapted} riemannian metric as one where the contact form $\alpha$ has length 1 and 
\bee \label{e:cah}
d\alpha=2\ast\alpha,
\ee
with $\ast$ denoting the hodge star operator associated with the metric. While this concept has been around somewhat longer in the more general form of a contact metric structure \cite{blair}, Chern and Hamilton showed that any oriented contact 3-manifold admits an adapted riemannian metric. We will adopt a slightly more general definition of adapted metrics by replacing (\ref{e:cah}) by
\[
d\alpha=2\lambda\ast\alpha,
\]
for any non-zero constant $\lambda$. Of course a rescaling of $\alpha$ and the metric can set $\lambda$ to 1, but it will be useful to leave $\lambda$ free for what follows.  We can think of this as a condition on the metric given a contact structure, or as a condition on the contact structure given the metric.

The purpose of this paper is to consider smooth contact structures adapted to Einstein metrics on $R^3$.  Our main result is:

\begin{thm}  \label{t:main} A hyperbolic metric does not admit an adapted contact structure.  Every smooth contact structure on $R^3$ adapted to a flat metric is contact isometric to $R^3$  with the standard metric and adapted contact 1-form. Every smooth contact structure on $R^3$ adapted to an elliptic metric is contact isometric to an open subset of $S^3$ with the round metric and standard adapted contact 1-form.
\end{thm}

Here, two smooth riemannian manifolds ($M_1$, $g_1$) and ($M_2$, $g_2$) with adapted contact structures given by contact 1-forms $\alpha_1$ and $\alpha_2$, respectively, are {\it contact isometric} if there exists a diffeomorphism $\phi:M_1\rightarrow M_2$ such that $\phi^*g_2=g_1$ and $\phi^*\alpha_2=\alpha_1$.

The main tool we will use is that of spin coefficients, first introduced by Newman and Penrose \cite{nap} in the lorentz setting. We adapt this technique, which was refined to the stationary case by Perjes \cite{perj}, to the 3-dimensional riemannian setting.  The basic idea is to replace the metric coefficients by the connection coefficients, taken in complex combinations motivated by the spinor formalism, as the primary variables.  In the presence of a preferred congruence of curves these spin coefficients, while breaking covariance, can be interpreted in terms of the acceleration, twist, divergence and shear of the congruence.  For a contact structure with contact 1-form there is a preferred vector field, called the Reeb vector field.  When the contact structure is adapted to a metric the Reeb vector field is orthogonal to the contact planes.  Moreover, we show that the integral curves of the Reeb vector field are characterized by certain geometric properties: 
 
\noindent {\bf Theorem 3.1}
{\it
Let $\xi$ be an contact structure with adapted metric $g$.  Then the congruence of curves orthogonal to $\xi$ are
\begin{description}
\item[(i)] geodesic
\item[(ii)] divergence-free
\item[(iii)] of constant non-zero twist.
\end{description}
Conversely, if a metric $g$ admits a geodesic congruence with constant (non-zero) twist, then there exists a contact structure adapted to $g$.  In addition, this congruence is divergence-free.
}

An adapted contact structure gives rise to a pseudohermitian structure \cite{tan} \cite{web}.  For such a structure one can define the Tanaka-Webster curvature and torsion.  In terms of the spin coefficients on a 3-manifold we show that the pseudohermitian torsion is the shear, while we identify the curvature as certain derivatives of the spin coefficients.  These are related to the scalar curvature of the riemannian metric in such a way that:

\noindent {\bf Proposition 3.2}
{\it
If any two of the scalar curvature, the Tanaka-Webster curvature and the torsion are constant, then all three are.
 }

While every contact manifold admits an adapted metric, we will look at contact structures for which the adapted metric is Einstein. In fact, it is known that a locally symmetric 3-manifold with adapted contact structure must be Einstein \cite{bas}.  In three dimensions being Einstein is equivalent to the metric being of constant curvature.  That is, the manifold is covered by $H^3$, $R^3$ or $S^3$ with the standard metrics. Thus we consider the local moduli space of contact structures adapted to these three geometries. We prove the following:

\noindent {\bf Theorem 4.1}
{\it
The torsion of a contact structure adapted to an Einstein 3-manifold is constant.  In addition, the manifold is elliptic if and only if the torsion is zero and is flat if and only if the torsion is non-zero.  A hyperbolic manifold does not admit an adapted contact structure.
}

This paper is arranged as follows.  In the next Section we introduce the Newman-Penrose-Perjes spin coefficient formalism tailored to the three dimensional riemannian setting.  In the following Section we apply this formalism to contact structures and identify the Tanaka-Webster curvature and torsion in terms of the spin coefficients.

In Section 4 we apply the technique to contact structures adapted to Einstein metrics and prove the main theorem. In the final Section we discuss further applications of this technique in contact geometry.

Throughout we use the abstract index notation and the summation convention, and all indices are raised and lowered by the riemannian metric $g$ on the 3-manifold $M$.

\section{The Newman-Penrose-Perjes Formalism }

Let ($M$, $g$) be a riemannian 3-manifold and consider an orthonormal frame $\{e_{0}=Z_{0}, e_1, e_2\}$, and dual basis of 1-forms $\{\theta^{0}, \theta^1, \theta^2\}$.  Introduce a complex frame $\{Z_{0}, Z_+, Z_-\}$ by 
\[
Z_+=\frac{1}{\sqrt{2}}\left(e_1-i\;e_2\right) \qquad\qquad Z_-=\frac{1}{\sqrt{2}}\left(e_1+i\;e_2\right),
\]
with dual basis of 1-forms $\{\theta^{0}, \theta^+, \theta^-\}$.

We define the complex {\it spin coefficients} by
\[
\gamma_{mnp}=\nabla_jZ_{mi}\;Z_n^iZ_p^j,
\]
where $\nabla$ is the covariant derivative associated with $g$ and the indices $m$, $n$, $p$ range over $0$, $+$, $-$.  Thus
\[
\gamma_{mnp}=-\gamma_{nmp}.
\]
Alternatively, the 1-forms and the spin coefficients are related by Cartan's first equation of structure
\bee \label{e:cartan1}
d\theta^m=\theta^n\wedge\omega_n^{\;\;m},
\ee
where
\[
\omega_{nm}=\gamma_{nmp}\theta^p.
\]
We break covariance and introduce the complex optical scalars
\[
\gamma_{+0-}=\rho \qquad\qquad \gamma_{+0+}=\sigma \qquad\qquad \gamma_{+--}=\tau
\]
\[
\gamma_{+00}=\kappa \qquad\qquad \gamma_{+-0}=\epsilon.
\]
Geometrically, $\kappa=0$ if and only if $e^i_{0}$ is tangent to a geodesic congruence in $M$. In order to interpret the other coefficients geometrically, consider a bundle of geodesics (often called a {\it geodesic congruence}) with affine parameter $r$ along the geodesics.  Consider a {\it connecting vector} $V^i$ along the geodesic $\gamma$, that is, $V^i$ joins points on $\gamma$ with those of the same parameter value on a neighboring geodesic $\gamma'$.  If the orthonormal frame vector $e^i_{0}$ is chosen to be tangent to the congruence, then $V^i$ is Lie transported along the rays:
\bee \label{e:liedrag}
{\cal{L}}_{e_{0}}V=0 \qquad\qquad \Longleftrightarrow \qquad\qquad e_{0}^j\nabla_j V^i=V^j\nabla_je_{0}^i.
\ee
Without loss of generality we can choose the parameter $r$ so that $e_{0}^iV_i=0$ at some point, and then it is not hard to show that this holds along the whole of $\gamma$.  Thus, in terms of the orthonormal frame $\{e_{0}, Z_+, Z_-\}$, we have
\[
V^i=\zeta Z_-^i+\overline{\zeta}Z^i_+,
\]
for some complex function $\zeta$.  Now equation (\ref{e:liedrag}) reduces to the single complex equation
\bee\label{e:optics}
\frac{d \zeta}{dr}=-\rho\zeta-\sigma\overline{\zeta}.
\ee
We use this equation to track the behavior of the geodesic $\gamma'$ relative to $\gamma$.  Let $\rho=\theta+\nu i$ and $\sigma=|\sigma|e^{2i\phi}$. If $\sigma=\nu=0$, then equation (\ref{e:optics}) reads
\[
\frac{d \zeta}{dr}=-\theta\zeta,
\]
so the real part of $\rho$ measures the rate of divergence or contraction of the geodesics.  Therefore $Re(\rho)$ is referred to as the {\it divergence} of the congruence.  On the other hand if $\sigma=\theta=0$,then
\[
\frac{d \zeta}{dr}=-i\nu\zeta.
\]
Thus the imaginary part of $\rho$ measures the rotation of neighboring geodesics, and so we call $Im(\rho)$ the {\it twist} of the congruence.  Finally, if $\rho=0$ the propagation equation is
\[
\frac{d \zeta}{dr}=-|\sigma|e^{2i\phi}\overline{\zeta},
\]
i.e. the derivative of $\zeta$ is a real multiple of $\zeta$ when $arg(\zeta)=\phi$, $\phi+\pi$ or $\phi\pm\hf\pi$.  In the first two cases we get contraction towards the origin, while in the last two cases we get dilation away from the origin - a circle lying in the plane orthogonal to the geodesic will propagate to an ellipse.  Thus $|\sigma|$ is a measure of the eccentricity of the ellipse and it is called the {\it shear} of the congruence, while $\phi$ measures the inclination of the ellipse in this plane.

It is not hard to show that, in the case where the congruence is geodesic, the divergence can be expressed by 
\[
Re(\rho)={\hf}\nabla_i e_0^i,
\]
the twist by
\[
Im(\rho)={\hf}\left[\nabla_{[i}e_{0\;j]}\nabla^i e_0^j\right]^{\shf},
\] 
and the shear by
\[
|\sigma|={\sdv{1}{\sqrt{2}}}\left[\nabla_{(i}e_{0\;j)}\nabla^i e_0^j -{\hf}\left(\nabla_i e_0^i \right)^2\right]^{\shf}.
\]

The Riemann tensor of $g$ is given in terms of the spin coefficients by
\[
R_{mnpq}=\gamma_{mnp,q}-\gamma_{mnq,p}+\gamma^r_{\;\;mq}\gamma_{rnp} - \gamma^r_{\;\;mp}\gamma_{rnq} + \gamma_{mnr}\left(\gamma^r_{\;\;pq} - \gamma^r_{\;\;qp}\right),
\]
where a coma subscript represents ordinary differentiation in the indicated direction.

We introduce the differential operators
\[
D=Z_{0}^i\;\frac{\partial}{\partial x^i} \qquad \qquad \delta=Z_+^i\;\frac{\partial}{\partial x^i}  \qquad \qquad \overline{\delta}=Z_-^i\;\frac{\partial}{\partial x^i},
\]
and if $f$ is any function, the commutators of these operators work out to be
\bee \label{e:comm1}
(D\delta-\delta D)f=[(\overline{\rho}+\epsilon)\delta +\sigma\overline{\delta}+\kappa D]f
\ee
\bee \label{e:comm2}
(\delta\overline{\delta}-\overline{\delta}\delta)f=[\overline{\tau}\overline{\delta} -\tau\delta+(\overline{\rho}-\rho) D]f.
\ee
The components of the Ricci tensor in terms of the spin coefficients are:
\bee \label{e:ric00}
R_{00}=D\rho+D\overline{\rho}-\overline{\delta}\kappa-\delta\overline{\kappa}+\tau\kappa+\overline{\tau}\;\overline{\kappa}-2\kappa\overline{\kappa}                       
        -2\sigma\overline{\sigma}-\rho^2-\overline{\rho}^2
\ee
\bee \label{e:ric++}
R_{++}=-\delta\kappa+D\sigma
-2\epsilon\sigma-\overline{\tau}\kappa-\kappa^2 
        -\sigma\overline{\rho}-\rho\sigma 
\ee
\bee \label{e:ric0+}
R_{0+}=-\overline{\delta}\sigma+\delta\rho+2\tau\sigma+\kappa\rho-\kappa\overline{\rho} 
\ee
\bee \label{e:ric0-}
R_{0-}=-\overline{\delta}\epsilon+D\tau +\kappa\overline{\sigma}-\rho\overline{\kappa}+\epsilon\tau
        -\epsilon\overline{\kappa}+\overline{\tau}\;\overline{\sigma}-\tau\rho 
\ee
\bee \label{e:ric+-}
R_{+-}= -\overline{\delta}\kappa+D\rho+\delta\tau+\overline{\delta}\;\overline{\tau}+\epsilon\rho
        -\epsilon\overline{\rho}-\kappa\overline{\kappa}+\kappa\tau-\rho\overline{\rho}-\rho^2
        -2\tau\overline{\tau},
\ee
while the scalar curvature is
\bee \label{e:rsc}
{\shf}R=-2\delta\overline{\kappa}+2D\overline{\rho}+\delta\tau+\overline{\delta}\overline{\tau}-2\kappa\overline{\kappa}+2\overline{\kappa}\;\overline{\tau}-2\overline{\rho}^2-\sigma\overline{\sigma}+\epsilon\rho-\epsilon\overline{\rho}-\rho\overline{\rho}-2\tau\overline{\tau} .
\ee
In addition, we have the following identities from the symmetries of the curvature tensor
\bee \label{e:id1}
D\rho-\overline{\delta}\kappa+\kappa\tau-\rho^2 = D\overline{\rho}-\delta\overline{\kappa}
        +\overline{\kappa}\;\overline{\tau}-\overline{\rho}^2
\ee
\bee \label{e:id2}
\delta\overline{\sigma}-\overline{\delta}\overline{\rho}-\overline{\tau}\;\overline{\sigma}-\overline{\kappa}\;\overline{\rho} = 
   \overline{\delta}\epsilon-D\tau -\kappa\overline{\sigma} -\epsilon\tau+\epsilon\overline{\kappa} +\tau\rho.
\ee

\section{Adapted Contact Structures}

Following (and slightly generalizing) Chern and Hamilton \cite{cah} we will say that a contact structure $\xi$ on a riemannian  3-manifold is {\it adapted} to the metric $g$ on $M$ if there exists a contact form $\alpha$ for $\xi$ such that
\bee \label{e:adapted}
d\alpha=2 \lambda\ast\alpha \qquad \qquad g(\alpha,\alpha)=1,
\ee
where $\ast$ is the hodge star operator associated with $g$ and $\lambda$ is a non-zero constant.   This condition has also been studied by Nicolaescu \cite{nic} , who referred to it as $\lambda$-adapted.  The sign of $\lambda$ can be changed by switching the orientation of the metric.  We will assume throughout that the orientation is chosen so that $\lambda$ is positive.

\noindent{\bf Examples}

\noindent {\bf (a)} The standard example of an adapted contact structure on $R^3$ is the flat metric and contact 1-form
\bee \label{e:constan}
\alpha=\sin(2\lambda z)dx+\cos(2\lambda z)dy.
\ee

\noindent {\bf (b)} If we consider the 3-sphere of radius $c$ in $S^3\subset R^4$ given by $x^2+y^2+z^2+w^2=c^2$, then the round metric and contact 1-form
\bee \label{e:3sphere}
\alpha=\frac{1}{c}(xdy-ydx+zdw-wdz),
\ee
form an adapted contact structure with $c=\lambda^{-1}$.

\begin{thm}
Let $\xi$ be an adapted contact structure.  Then the congruence of curves orthogonal to $\xi$ are
\begin{description}
\item[(i)] geodesic
\item[(ii)] divergence-free
\item[(iii)] of constant non-zero twist.
\end{description}

Conversely, if $M$ has a metric $g$ which admits a geodesic congruence with constant (non-zero) twist, then there exists a contact structure adapted to $g$.  In addition, this congruence is divergence-free.

\end{thm}

\noindent{\bf Proof:}:

Equation (\ref{e:adapted}) means that
\beer
\nabla_i\alpha_j - \nabla_j\alpha_i &=&2 \lambda \epsilon_{ij}^{\;\;\;\; k} \alpha_k \nn\\
\Rightarrow\quad \alpha^i \nabla_i\alpha_j - \alpha^i \nabla_j\alpha_i &=& 0 \nn\\
\Rightarrow\quad \alpha^i \nabla_i\alpha_j &=& \hf\left(|\alpha|_g^2\right)_{,j} \nn\\
&=& 0, \nn
\eer
so the congruence is geodesic (and affinely parametrised), proving part (i).

To prove (ii) take the exterior derivative of the first of equation (\ref{e:adapted})
\[
0=dd\alpha= 2\lambda d*\alpha \qquad\Rightarrow \qquad *d*\alpha=0,
\]
so that $\alpha$ is divergence-free.

Finally,
\beer
Im(\rho)&=&\hf \left[\nabla_{[i}\alpha _{j]}\nabla^i \alpha ^j\right] ^{\shf}\nn\\
&=& \hf\left(\hf \nabla_{[i}\alpha_{j]}\nabla^{[i}\alpha^{j]}\right)^{\shf}\nn\\
&=&\lambda.|*\alpha|_g \nn\\
&=&\lambda,\nn
\eer
so the twist is constant. 

Conversely, suppose that $\alpha^i$ is the unit 1-form dual to the vector field tangent to the geodesic congruence.  Then
\[
\alpha^j\nabla_{[i}\alpha_{j]} = \hf\left[|\alpha|_g^2\right]_{,j}-\alpha^j\nabla_j\alpha_i=0.
\]
Thus $\alpha^i$ is in the null space of the 2-form $d\alpha$.  Form an orthonormal basis of 1-forms $\{\alpha, \theta^1,\theta^2\}$.  Then the space of 2-forms has a basis of the form $\{\alpha\wedge\theta^1, \alpha\wedge\theta^2, \theta^1\wedge\theta^2\}$ and we can resolve $d\alpha$ as
\[
d\alpha=f\alpha\wedge\theta^1+g \alpha\wedge\theta^2+h \theta^1\wedge\theta^2,
\]
for some functions $f$, $g$ and $h$.  Now, since $\alpha^i$ is in the null space of  $d\alpha$, $f=g=0$ and we must have
\[
d\alpha=h \theta^1\wedge\theta^2=h*\alpha.
\]
Finally, since the twist is constant
\[
Im(\rho)={\hf}\left(|d\alpha|_g^2\right)^{\shf}=\frac{|h|}{2}=constant.
\]
Thus, $|h|=2\lambda$ and we have an adapted contact structure.

Alternatively, we can show that a geodesic congruence with constant (non-zero) twist is divergence-free as follows.  Consider the identity (\ref{e:id1}) for a geodesic congruence ($\kappa=0$):
\[
D\rho -\rho^2 = D\overline{\rho}-\overline{\rho}^2.
\]
If we let $\rho=\theta+\nu i$, with constant $\nu$ then this says that $\theta.\nu=0$. Since $\nu\neq0$ we conclude that $\theta=0$, i.e. the congruence is divergence-free. 

$\qquad$ {\it q.e.d.}

\vspace{0.2in}

A {\it pseudohermitian structure} on $M$ is a triple ($\xi, J, \alpha$) such that $\xi$ is a contact structure with contact form $\alpha$ and $J:\xi\rightarrow \xi$ is a smooth endomorphism such that $J^2=-Id_\xi$.  That is, $J$ is a CR-structure on $\xi$.  A pseudohermitian structure gives rise to a riemannian metric on $M$ as follows.

Choose a complex vector field $Z_+$ which is an eigenvector of $J$ with eigenvalue $i$, and let $Z_-=\overline{Z_+}$.  We can complete these to  form a frame $\{Z_{0}, Z_+, Z_-\}$ by adding the {\it Reeb vector field} uniquely defined by the properties
\[
\alpha(Z_{0})=1 \qquad d\alpha(Z_{0}, \cdot)=0.
\]

We denote the dual basis of 1-forms by $\{\theta^{0}, \theta^+, \theta^-\}$.  Moreover, $Z_+$ can be chosen so that
\bee \label{e:pseudo}
d\theta^{0}= 2\lambda i \theta^+\wedge\theta^-,
\ee
for some constant $\lambda$.  On the other hand
\bee \label{e:torsion}
d\theta^+=\theta^+\wedge\omega+A\theta^{0}\wedge\theta^-,
\ee
where $A$ is the {\it pseudohermitian torsion}.

We can associate a riemannian metric with the pseudohermitian structure by splitting $Z_+$ and $Z_-$ into real and imaginary parts:
\[
Z_+=\frac{1}{\sqrt{2}}\left(e_1-i\;e_2\right) \qquad\qquad Z_-=\frac{1}{\sqrt{2}}\left(e_1+i\;e_2\right),
\]
and defining an orthonormal frame by $\{e_{0}=Z_{0}, e_1, e_2\}$.  That is, the metric is defined by
\[
g^{ij}= e_{0}^i\; e_{0}^j + e_1^i\; e_{0}^j+ e_2^i\; e_2^j,
\]
while the dual basis of 1-forms will be denoted by $\{\theta^{0}, \theta^1, \theta^2\}$.  

The {\it Tanaka-Webster curvature} $W$ of the pseudohermitian structure is defined by
\bee \label{e:webcurv}
d\omega=2\lambda W\theta^+\wedge\theta^-+2i Im(\nabla_{Z_-}\overline{A}\;\theta^+\wedge\theta^{0}).
\ee

The real version of equation (\ref{e:pseudo}) is simply the adapted condition
\[
d\alpha=2\lambda \ast\alpha.
\]

This construction is reversible, so that a pseudohermitian structure on a 3-manifold is equivalent to a contact structure adapted to a riemannian metric.

We now identify the torsion and Tanaka-Webster curvature  of a pseudohermitian structure in terms of the spin coefficients.  To do this, consider equation (\ref{e:cartan1}) with $m=+$
\beer
d\theta^+&=& \theta^n\wedge\omega_n^{\;\;+}\nn\\
&=& \theta^+\wedge\omega_+^{\;\;+}+\theta^{0}\wedge\left(\gamma_{0-0}\theta^{0} +\gamma_{0--}\theta^-+\gamma_{0-+}\theta^+\right)\nn\\
&=& \theta^+\wedge\left(\omega_{+-}-\gamma_{0-+}\theta^{0}\right) +\gamma_{0--}\theta^{0}\wedge\theta^-\nn\\
&=& \theta^+\wedge\left(\omega_{+-}+\overline{\rho}\theta^{0}\right) -\overline{\sigma}\theta^{0}\wedge\theta^-.\nn
\eer
Comparing this with equation (\ref{e:torsion}) we see that the pseudohermitian torsion is determined by the shear:
\[
A=-\overline{\sigma},
\]
and 
\beer
\omega&=&\omega_{+-}+\overline{\rho}\theta^{0}\nn\\
&=&\left(\epsilon+\overline{\rho}\right)\theta^{0}-\overline{\tau}\theta^++\tau\theta^-.\nn
\eer
A calculation shows that
\[
d\omega=\left(\delta\tau+\overline{\delta}\overline{\tau}-2\overline{\tau}\tau -2\lambda\rho i\right)\theta^+\wedge\theta^-+2i\;Im\left(D\overline{\tau}+\delta\epsilon-\overline{\tau}(\epsilon+\overline{\rho})+\tau\sigma\right)\theta^+\wedge\theta^{0}.
\]
Comparing this with equation (\ref{e:webcurv}) the Tanaka-Webster curvature is given by
\[
2\lambda W=\delta\tau+\overline{\delta}\overline{\tau}-2\overline{\tau}\tau +2\lambda^2.
\]
Note that we also have the identity
\[
D\overline{\tau}+\delta\epsilon-\overline{\tau}(\epsilon+\overline{\rho})-\tau\sigma=-\overline{\delta}\sigma,
\]
which is just equation (\ref{e:id2}) for an adapted contact structure.

Thus we have that the Tanaka-Webster curvature and the scalar curvature of the associated metric are related by
\[
\hf R= 2\lambda W-\lambda^2 -|\sigma|^2.
\]

\begin{prop}
If any two of the scalar curvature, the Tanaka-Webster curvature and the torsion are constant, then all three are. In addition, if the scalar curvature $R$ is positive, then so is the Tanaka-Webster curvature $W$.
\end{prop}

\section{Einstein 3-Manifolds}

In this section we look at the question of whether there exists contact structures adapted to Einstein 3-manifolds:
\[
R_{ij}={\sdv{1}{3}}Rg_{ij}.
\]
Since the Ricci curvature determines the full Riemann curvature in 3 dimensions this is equivalent to asking whether there exists adapted contact structures on 3-manifolds of constant curvature.

Our considerations will be local in this section.  We will use the Newman-Penrose-Perjes formalism to completely integrate the field equations.  We have seen how an adapted contact  structure gives rise to a divergence-free, constantly twisting congruence of geodesics. That is
\[
\kappa=0 \qquad\qquad \rho=- \overline{\rho}=\lambda i=constant.
\]
In addition, we can choose the frame so that $\epsilon=0$ and still have the remaining freedom to transform
\bee \label{e:triad}
Z_+\rightarrow e^{iC_{0}}Z_+ \qquad\qquad Z_-\rightarrow e^{-iC_{0}}Z_-,
\ee 
where $C_{0}$ is a real function which is independent of the affine parameter along the geodesics.  These simplifications mean that equations (\ref{e:ric00}) to (\ref{e:rsc}) become
\bee \label{e:ric00pha}
R_{00}=-2\sigma\overline{\sigma}-2\rho^2
\ee
\bee \label{e:ric++pha}
R_{++}= D\sigma
\ee
\bee \label{e:ric0+pha}
R_{0+}=-\overline{\delta}\sigma +2\tau\sigma  
\ee
\bee \label{e:ric0-pha}
R_{0-}= D\tau +\overline{\tau}\;\overline{\sigma}-\tau\rho 
\ee
\bee \label{e:ric+-pha}
R_{+-}= \delta\tau+\overline{\delta}\;\overline{\tau}  -2\tau\overline{\tau}
\ee
\bee \label{e:rscpha}
{\shf}R= \delta\tau+\overline{\delta}\overline{\tau}-\rho^2-\sigma\overline{\sigma}-2\tau\overline{\tau} .
\ee
The first of the identities (\ref{e:id1}) is identically satisfied while the second becomes
\bee \label{e:id2ph}
\delta\overline{\sigma} -\overline{\tau}\;\overline{\sigma} =-D\tau +\tau\rho.
\ee
For an Einstein manifold we have that
\[
R_{++}=R_{0-}=R_{0+}=0 \qquad R_{00}=R_{+-}={\sdv{1}{3}}R \qquad R=constant.
\]
Substituting this in (\ref{e:ric00pha}) to (\ref{e:id2ph}) we get that 
\bee \label{e:ric00ph}
\sigma\overline{\sigma}+\rho^2=-{\sdv{1}{6}}R
\ee
\bee \label{e:ric++ph}
D\sigma=0
\ee
\bee \label{e:ric0+ph}
\overline{\delta}\sigma -2\tau\sigma=0
\ee
\bee \label{e:ric0-ph}
-D\tau -\overline{\tau}\;\overline{\sigma}+\tau\rho =0
\ee
\bee \label{e:ric+-ph}
-\delta\tau-\overline{\delta}\;\overline{\tau}+2\tau\overline{\tau}=-{\sdv{1}{3}}R.
\ee
Equation (\ref{e:ric++ph}) tells us that the shear (or torsion) $\sigma$ is constant along the geodesics.  Thus, by a transformation of the type (\ref{e:triad}) we can make $\sigma$ real, and by (\ref{e:ric00ph}), it will be constant.  From (\ref{e:ric0+ph}) then, we must have either $\sigma=0$ or $\tau=0$.  The latter, by equation (\ref{e:ric+-ph}), can only happen if the metric is flat ($R=0$), while according to equation (\ref{e:ric00ph}) the former can only happen when the metric is elliptic ($R>0$).  Thus we have the following:

\begin{thm} \label{t:einstein}
The torsion of a contact structure adapted to an Einstein 3-manifold is constant:
\[
|\sigma|=\left|\sqrt{ \lambda^2- {\sdv{1}{6}}R }\right|.
\]
In addition, the manifold ($M$, $g$) is elliptic if and only if the torsion is zero with $R=6\lambda^2$, and is flat if and only if the torsion is non-zero with $|\sigma|=\lambda$.  A hyperbolic manifold does not admit an adapted contact structure.
\end{thm}

\begin{cor} 
The Tanaka-Webster curvature of a contact structure adapted to an Einstein 3-manifold is constant:
\[
W= {\dv{1}{3\lambda}}R+\lambda.
\]
\end{cor}

Choose coordinates ($x^{0}=r,x^a$) where $r$ is an affine parameter along the geodesics and $a=1,2$.  Then
\[
D=\frac{\partial}{\partial r} \qquad \delta=\Omega\frac{\partial}{\partial r}+\eta^a\frac{\partial}{\partial x^a},
\]
for some complex functions $\Omega$ and $\eta^a$.  The remaining coordinate freedom we have is shifting the origin
\bee \label{e:coord1}
r'=r+r_{0}(x^a),
\ee
and relabelling the geodesics
\bee \label{e:coord2}
x^{a'}=x^{a'}(x^b).
\ee
The commutation relations (\ref{e:comm1}) and (\ref{e:comm2}) yield
\[
D\Omega=-\rho\Omega+\sigma\overline{\Omega}
\]
\[
D\eta^a=-\rho\eta^a+\sigma\overline{\eta}^a
\]
\[
\delta\overline{\Omega}-\overline{\delta}\Omega=\overline{\tau}\;\overline{\Omega}-\tau\Omega-2\rho
\]
\[
\delta\overline{\eta}^a-\overline{\delta}\eta^a=\overline{\tau}\;\overline{\eta}^a-\tau\eta^a.
\]
We will now deal with the flat and elliptic cases separately.

\subsection{Flat 3-Manifolds}

When $R=0$ we have $\sigma=\lambda=|\rho|$, $\tau=0$ and the equations remaining to be solved are
\bee \label{e:comm1fta}
D\Omega=-i\lambda \Omega+\lambda \overline{\Omega}
\ee
\bee \label{e:comm1ftb}
D\eta^a=-i\lambda \eta^a+\lambda \overline{\eta}^a
\ee
\bee \label{e:comm2fta}
\delta\overline{\Omega}-\overline{\delta}\Omega= -2\lambda i
\ee
\bee \label{e:comm2ftb}
\delta\overline{\eta}^a-\overline{\delta}\eta^a=0.
\ee
Applying $D$ to the first two of these we find that
\[
DD\Omega=0 \qquad DD\eta^a=0.
\]
Thus $\Omega=\alpha_{0}+\beta_{0}r$ where $\alpha_{0}$ and $\beta_{0}$ are complex functions of $x^a$ only.  From here on a subscript 0 will mean that the scalar is independent of $r$.  Putting this back into (\ref{e:comm1fta}) we find that
\[
\alpha_{0}=a_{0}+(b_{0}\lambda^{-1}-a_{0})i \qquad\qquad \beta_{0}=b_{0}(1-i),
\]
where $a_{0}$ and $b_{0}$ are real functions.  Similarly we find from equation (\ref{e:comm1ftb}) that
\[
\eta^a=A_{0}^a+(B_{0}^a\lambda^{-1}-A_{0}^a)i+B_{0}^a(1-i)r,
\]
for some real functions $A^a_{0}$ and $B_{0}^a$.  We can view these as the components of vector fields
\[
\vec{A}=A_{0}^a\frac{\partial}{\partial x^a} \qquad\qquad \vec{B}=B_{0}^a\frac{\partial}{\partial x^a}.
\]
The remaining two equations (\ref{e:comm2fta}) and (\ref{e:comm2ftb}) now reduce to
\[
\vec{A}(b_{0}) - \vec{B}(a_{0})=\lambda^2+b_{0}^2
\]
\[
\left[\vec{A},\vec{B}\right]=b_{0}\vec{B}.
\]
We still have the freedom of making the change of coordinates (\ref{e:coord1}) and (\ref{e:coord2}).  In order to utilize these we need to treat the cases $b_{0}=0$ and $b_{0}\neq 0$ separately.

\subsubsection{$b_{0}=0$}

In this case $\Omega=a_{0}(1-i)$ and the remaining equations are
\[
\vec{B}(a_{0})=-\lambda^2
\]
\[
\left[\vec{A},\vec{B}\right]=0.
\]
If $\vec{A}\propto \vec{B}$ then the metric turns out to be degenerate and so we will assume that this is not the case.  Thus, by the second equation we can choose local coordinates ($u, v$) so that
\[
\vec{A}=\frac{\partial}{\partial u} \qquad\qquad \vec{B}=\frac{\partial}{\partial v},
\]
and then the first equation  says that 
\bee \label{e:finalflata}
a_{0}=-\lambda^2v+f,
\ee
for some function $f=f(u)$.

We can assemble the local form for the metric:
\[
g_{ij}=\left[\begin{array} {ccc}
1 & -a_{0} & 0 \nn\\
-a_{0} & {\shf}(2\lambda ^2r^2-2\lambda r+1) +a_{0}^2& -{\shf}(2\lambda ^2r-\lambda ) \nn\\
0 &  -{\shf}(2\lambda ^2r-\lambda ) & \lambda ^2 \nn\\
\end{array}\right],
\]
where $a_{0}$ is given by (\ref{e:finalflata}).

This is contact  isometric to the standard adapted contact structure on $R^3$.  We can see this by considering the diffeomorphism $\phi(r,u,v)= (x,y,z)$ defined by
\[
x=\left(r-\frac{1}{2\lambda}\right)\sin(\lambda u)-\left(\lambda v-\frac{f}{\lambda}\right)\cos(\lambda u)-\frac{1}{\lambda}\int \cos(\lambda u)\frac{df}{du}du
\]
\[
y=\left(r-\frac{1}{2\lambda}\right)\cos(\lambda u)+\left(\lambda v-\frac{f}{\lambda}\right)\sin(\lambda u)+\frac{1}{\lambda}\int \sin(\lambda u)\frac{df}{du}du
\]
\[
z=\frac{u}{2}.
\]
The metric pulls back to $(\phi^{-1})^*g_{ij}=\delta_{ij}$ and the contact 1-form pulls back to the standard form
\[
(\phi^{-1})^*\alpha=\sin(2\lambda z)dx+\cos(2\lambda z)dy.
\]

\subsubsection{$b_{0}\neq 0$}
In this case we can utilize the change of coordinates (\ref{e:coord1}) to set $a_{0}=0$.  Then by the remaining freedom (\ref{e:coord2}) we can set
\[
\vec{A}=f\frac{\partial}{\partial u} \qquad\qquad \vec{B}=g\frac{\partial}{\partial v},
\]
for coordinates ($r, u, v$), where $f$ and $g$ are some functions of $u$ and $v$.  The equations to be solved are then
\[
f\frac{\partial b_{0}}{\partial u}=\lambda^2+b_{0}^2
\]
\[
f\frac{\partial g}{\partial u}=gb_{0}
\]
\[
-g\frac{\partial f}{\partial v}=0.
\]
The last of these implies that $f=f(u)$ and a further coordinate change $u\rightarrow u'=u'(u)$ allows us to set $f=1$.  Then we get
\[
\frac{\partial b_{0}}{\partial u}=\lambda^2+b_{0}^2
\]
\[
b_{0}=\frac{1}{g}\frac{\partial g}{\partial u}.
\]
The solution to these equations is
\bee \label{e:finalflat1b}
\frac{1}{g}= D(v)\sin(\lambda u)+E(v)\cos(\lambda u)
\ee
\bee \label{e:finalflat2b}
b_{0}=-\lambda \;\frac{D(v)\cos(\lambda u)-E(v)\sin(\lambda u)}{D(v)\sin(\lambda u)+E(v)\cos(\lambda u)},
\ee
where $D$ and $E$ are arbitrary functions of $v$.  The metric turns out to be
\[
g_{ij}=\left[\begin{array} {ccc}
1 & 0 & -\frac{b_{0}}{g} \nn\\
0 & {\shf}(2\lambda ^2r^2-2\lambda r+1) & - \frac{1}{2g} (2\lambda ^2r-\lambda ) \nn\\
-\frac{b_{0}}{g} &  -\frac {1}{2g} (2\lambda ^2r-\lambda ) & \frac{1}{g^2}(\lambda ^2+b_{0}^2) \nn\\
\end{array}\right],
\]
where $g$ and $b_{0}$ are given by (\ref{e:finalflat1b}) and (\ref{e:finalflat2b}).

Again, pulling the metric and contact 1-form by a diffeomorphism $\phi$, in this case given by
\[
x=\left(r-\frac{1}{2\lambda}\right)\sin(\lambda u)- \lambda\int E(v)dv
\]
\[
y=\left(r-\frac{1}{2\lambda}\right)\cos(\lambda u)+\lambda\int D(v)dv
\]
\[
z=\frac{u}{2},
\]
the metric becomes $(\phi^{-1})^*g_{ij}=\delta_{ij}$ and the contact 1-form becomes 

\[
(\phi^{-1})^*\alpha=\sin(2\lambda z)dx+\cos(2\lambda z)dy.
\]

Thus we have the following:

\begin{thm} Every contact structure adapted on $R^3$  to a flat metric is contact isometric to the standard metric and contact 1-form
\[
\alpha=\sin(2\lambda z)dx+\cos(2\lambda z)dy.
\]
\end{thm}

This contact form is invariant under deck transformations $(x,y,z)\rightarrow (x+2\pi, y+2\pi, z+2\pi )$.  Thus we have

\begin{cor} Every contact structure adapted to a flat metric on $T^3$ is contact isometric to the standard metric and contact 1-form.
\end{cor}

\subsection{Elliptic 3-Manifolds}
We consider now the case where $R>0$.  As we saw in Theorem \ref{t:einstein} the pseudohermitian torsion vanishes $\sigma=0$, and we have $R=6|\rho|^2=6\lambda ^2$.

The equations then to be solved are

\bee \label{e:ell1}
D\tau = i\lambda\tau 
\ee
\bee \label{e:ell2}
\delta\tau+\overline{\delta}\;\overline{\tau}=2\tau\overline{\tau}+2\lambda ^2
\ee
\bee \label{e:comm1ella}
D\Omega=-i\lambda \Omega
\ee
\bee \label{e:comm1ellb}
D\eta^a=-i\lambda \eta^a
\ee
\bee \label{e:comm2ella}
\delta\overline{\Omega}-\overline{\delta}\Omega=\overline{\tau}\;\overline{\Omega}-\tau\Omega -2\lambda i
\ee
\bee \label{e:comm2ellb}
\delta\overline{\eta}^a-\overline{\delta}\eta^a=\overline{\tau}\;\overline{\eta}^a-\tau\eta^a.
\ee
The solutions to equations (\ref{e:ell1}), (\ref{e:comm1ella}) and (\ref{e:comm1ellb}) are
\[
\tau=\tau_{0}\;e^{i\lambda r}
\]
\[
\Omega=\Omega_{0}\;e^{-i\lambda r}
\]
\[
\eta^a=\eta_{0}^a\;e^{-i\lambda r}.
\]
We can choose coordinates ($x^2,x^3$) so that 
\[
\eta^a_{0}=P_{0}(\delta^a_2+i\delta^a_3),
\]
and a frame so that $P_{0}$ is a real non-zero function of ($x^2,x^3$).  The remaining equations (\ref{e:ell2}), (\ref{e:comm2ella}) and (\ref{e:comm2ellb}) become
\bee\label{e:finell1}
2P_{0}\left[ \frac{\partial \tau_{0}}{\partial\overline{z}}+  \frac{\partial \overline{\tau}_{0}}{\partial z}\right] +(\Omega_{0}\tau_{0}-\overline{\Omega}_{0}\;\overline{\tau}_{0})\lambda i=2\tau_{0}\overline{\tau}_{0}+2\lambda ^2
\ee
\bee\label{e:finell2}
2P_{0}\left[ \frac{\partial \overline{\Omega}_{0}}{\partial\overline{z}} -  \frac{\partial \Omega_{0}}{\partial z}\right] + 2\Omega_{0}\;\overline{\Omega}_{0} \lambda i =\overline{\tau}_{0}\;\overline{\Omega}_{0}-\tau_{0}\Omega_{0}-2\lambda i
\ee
\bee\label{e:finell3}
2\frac{\partial P_{0}}{\partial \overline{z}}=\overline{\tau}_{0}-\Omega_{0}\lambda i,
\ee
where we have introduced the complex coordinates $z=x^2+ix^3$ and $\overline{z}=x^2-ix^3$ and so
\[
\frac{\partial}{\partial z}=\frac{1}{2}\left(\frac{\partial}{\partial x^2}-i\frac{\partial}{\partial x^3}\right)
\qquad\qquad
\frac{\partial}{\partial \overline{z}}=\frac{1}{2}\left(\frac{\partial}{\partial x^2}+ i\frac{\partial}{\partial x^3}\right).
\]
Now, differentiating (\ref{e:finell3}) and using (\ref{e:finell1}) and (\ref{e:finell2}) yields
\[
2P_{0}^2\frac{\partial^2 \ln P_{0}}{\partial z \partial \overline{z}}=\lambda ^2.
\]
The solution to this equation (using a change of complex coordinates ($z$, $\overline{z}$)) is
\[
P_0=\frac{\lambda }{\sqrt{2}}(1+z\overline{z}).
\]

A relabelling of the rays $r\rightarrow r'=r+f(x^2,x^3)$ allows us to set the imaginary part of $\Omega_{0}$ to zero.  Then if ($x^2=u$, $x^3=v$), the final equation to be solved (equation (\ref{e:finell2})) is
\[
{\hf}(1+u^2+v^2)\frac{\partial \Omega_0}{\partial v}-v\Omega_0=-\frac{1}{\sqrt{2}},
\]
which has solution
\bee \label{e:om0}
\Omega_{0}=-\frac{1}{\sqrt{2}} \left[\frac{v}{1+u^2}+\frac{(1+u^2+v^2)}{(1+u^2)^{3/2}}
    \tan^{-1}\left(\frac{v}{\sqrt{1+u^2}}\right) 
      +f(u)(1+u^2+v^2)\right],
\ee
for an arbitrary function $f(u)$.

Finally we have the local form for the metric:
\[
g_{00}=1
\]
\[ 
g_{01}= g_{10}=-\frac{\sqrt{2}\Omega_{0}}{\lambda (1+u^2+v^2)}
\]
\[
g_{11}=\frac{2\Omega_{0}^2+1}{\lambda ^2(1+u^2+v^2)^2}
\]
\[
g_{22}=\frac{1}{\lambda ^2(1+u^2+v^2)^2},
\]
where $\Omega_0$ is given by (\ref{e:om0}).

This is contact isometric the standard adapted contact structure on $S^3$. To see this, first utilize the transformation $(r,u,v)\rightarrow (\tilde{\rho},\tilde{\theta},\tilde{\phi})$ given by:

\[
r+{\frac{1}{\lambda}}\int f(u) du=\frac{1}{\lambda}\left[ \tilde{\rho}-\frac{\cos\tilde{\phi}\tan\tilde{\theta}}{\sqrt{1+\cos^2\tilde{\phi}\tan^2\tilde{\theta}}}\tan^{-1}\left(\frac{\sin\tilde{\phi}\tan\tilde{\theta}}{\sqrt{1+\cos^2\tilde{\phi}\tan^2\tilde{\theta}}}\right)\right] 
\]
\[
u=\cos\tilde{\phi}\tan\tilde{\theta}
\]
\[
v=\sin\tilde{\phi}\tan\tilde{\theta}.
\]
This transforms the metric and 1-form to
\[
ds^2=\frac{1}{\lambda^2}\left[   d\tilde{\rho}^2+d\tilde{\theta}^2+\sin^2\tilde{\theta}(d\tilde{\phi}^2-2d\tilde{\phi}d\tilde{\rho})   \right] 
\]
\[
\alpha=\frac{1}{\lambda}\left[ d\tilde{\rho}-\sin^2\tilde{\theta}d\tilde{\phi}\right] .
\]
Finally, a transformation $(\tilde{\rho},\tilde{\theta},\tilde{\phi})\rightarrow (\rho,\theta,\phi)$ given by
\[
\tilde{\rho}=\tan^{-1}\left[\frac{\cos\rho}{\sin\rho\cos\theta}\right]
\]
\[
\tilde{\theta}=\cos^{-1}\left[\sqrt{\cos^2\rho+\sin^2\rho\cos^2\theta}\right]
\]
\[
\tilde{\phi}=\tan^{-1}\left[\frac{\cos\rho}{\sin\rho\cos\theta}\right]-\phi,
\]
produces the round metric of radius $\lambda^{-1}$ in standard co-ordinates
\[
ds^2=\frac{1}{\lambda^2}\left[ d\rho^2+\sin^2\rho(d\theta^2+\sin^2\theta d\phi^2)\right] ,
\]
and contact  structure 
\[
\alpha=\frac{1}{\lambda}\left[ -\cos\theta d\rho+\cos\rho\sin\rho\sin\theta d\theta-\sin^2\rho\sin^2\theta d\phi\right] ,
\]
which is precisely the contact structure given by (\ref{e:3sphere}). This completes Theorem \ref{t:main}.

\section{Conclusions}

We have seen how the Newman-Penrose-Perjes spin coefficient formalism can be applied to contact structures on 3-manifolds.  In the case where the contact structure is adapted to a riemannian metric with constant curvature this method is powerful enough to give a complete local description of both metric and contact structure.  

This method can also be applied to contact structures adapted to non-Einstein metrics.  For example, it can be shown that a Sasakian structure on a 3-manifold \cite{blair} is simply a contact structure adapted to a riemannian metric for which the shear of the geodesic congruence vanishes.  This simplification allows one to integrate the curvature conditions and uncover important properties of such manifolds.  This result, which we will report on in a future paper, serves to compliment the recent classification of closed Sasakian 3-manifolds by Geiges \cite{geiges} and their CR structures by Belgun \cite{belgun}.

Eliashberg \cite{elias1} has introduced the notion of tight and overtwisted contact structures.  While the classification of overtwisted contact structures has been completed, that of tight contact structures remains open. Theorem \ref{t:main} means that contact structures adapted to space forms (when they exist) are tight.  One question then is, what geometric properties have metrics adapted to overtwisted contact structures?

A further application of this technique would be to investigate contact structures adapted to the other 5 locally homogenous three dimensional geometries \cite{scott}.  Again it would be expected that the explicit forms of the local geometry should be obtainable from this method.

\end{document}